\documentclass[12pt]{article}

\usepackage{amsmath,amssymb}
\usepackage{verbatim}
\usepackage{color}
\usepackage{graphicx}
\usepackage{amsthm}
\usepackage{url}

\newtheorem{theorem}{Theorem}[section]
\newtheorem{thm}[theorem]{Theorem}
\newtheorem{corollary}[theorem]{Corollary}

\newtheorem{lemma}[theorem]{Lemma}

\newtheorem{proposition}[theorem]{Proposition}

\newtheorem{definition}[theorem]{Definition}

\newtheorem{remark}[theorem]{Remark}

\numberwithin{equation}{section}

\def\be{\begin{equation}}
\def\ee{\end{equation}}
\def\bes{\begin{equation*}}
\def\ees{\end{equation*}}

  \def\sC {{\mathcal C}}
 \def\sE {{\mathcal E}} \def\sF {{\mathcal F}}

\def\sP {{\mathcal P}}  
\def\sS {{\mathcal S}}  
\def\sV {{\mathcal V}}  \def\sX {{\mathcal X}}
\def\sY {{\mathcal Y}}

 \def\bE {{\mathbb E}}

 \def\bN {{\mathbb N}} 
\def\bP {{\mathbb P}}

 \def\bZ {{\mathbb Z}}

\def\sms{\smallskip}

\def\sm{\smallskip\noindent}

\def\ignore#1{}

\def\al {\alpha}
\def\lam {\lambda} \def\Lam {\Lambda} 
\def\eps{\varepsilon}
\def\th{\theta} 

\def\Gam{\Gamma}

\def\to {\rightarrow}

\def\pd {\partial}
\def\q{\quad} 
\def\dint{\int\kern-.6em\int}

\def \half {{\tfrac12}}

\def\wt{\widetilde}

\def\be{\begin{equation}}
\def\ee{\end{equation}}
\def\bes{\begin{equation*}}
\def\ees{\end{equation*}}
\def\ba{\begin{align}}
\def\ea{\end{align}}
\def\xxea{\end{align}}
\def\bas{\begin{align*}}
\def\eas{\end{align*}}

\def\proof{{\smallskip\noindent {\em Proof. }}}
\def\qed{{\hfill $\square$ \bigskip}}

\definecolor{dgreen}{rgb}{0, 0.6, 0.1}
\definecolor{dblue}{rgb}{0, 0.0, 0.6}
\definecolor{vdblue}{rgb}{0,.08, 0.45}
\definecolor{dred}{rgb}{0.7, 0.0, 0.0}
\definecolor{vdblue}{rgb}{0,.08, 0.45}

\definecolor{purple}{rgb}{0.6, 0.0, 0.6}
\definecolor{mytext}{rgb}{0.1, 0.1, 0.1}

\begin{document}

\font\titlefont=cmbx14 scaled\magstep1
\title{\titlefont  \vspace{-5ex} A branching process with contact tracing}

\author{
M. T. Barlow\footnote{Research partially supported by NSERC (Canada)}
}

\maketitle

%\ver

\begin{abstract}
We consider a supercritical branching process and define a contact tracing mechanism on 
its genealogical tree. We calculate the growth rate of the post tracing process, and give
conditions under which the tracing is strong enough to drive the process to extinction.


\end{abstract}

\section{Introduction} \label{sec:intro}

In this paper we present a simple model for contact tracing during  an epidemic. 
The epidemic is taken to be a standard (discrete time) Bienaym\'e-Galton-Watson branching process
$(Z_n, n \ge 0)$ 
with mean $\lam$; the case of interest is when the process is
supercritical, so $\lam>1$. (See \cite{AN,H} for background on these processes.)
We assume that $b$ generations after infection, an infected individual is detected with
probability $p$. If an individual is detected as infected, an attempt is made to trace all
this individual's contacts, both forwards and backwards.
 We assume that for each infection link the probability of a successful
trace is $\al$, and that it can be determined with probability 1
whether or not a traced individual has been infected. If a traced individual is detected as infected,
then in turn all contacts of that individual are traced, and this process is repeated throughout
the genealogy of the epidemic. It is assumed that detected individuals are quarantined
and so can be removed from the pool of infectious individuals. 
(See the next section for a more precise definition.)
We write $Z^{CT}$ for the branching process after contact tracing.

If $\al=1$ then all traces are successful, and as soon as one individual is detected 
every infected individual will be traced and isolated, so bringing the epidemic to an end.
For smaller values of $\al$ there is the possibility that the epidemic will remain supercritical:
we are interested in how large $\al$ needs to be to control the epidemic for a fixed $b$ and $p$.

Our main result is as follows. Let
\be \label{e:Gdef}
 G(u) = \sum_{k=0}^\infty p_k u^k
\ee
be the p.g.f. of the offspring distribution. Define sequences $(g_n)$, $(h_n)$ by
$g_0=h_0=1-p$, and
\begin{align}
   g_n(s) &= (1-p) G\big( 1-\al + \al g_{n-1}\big), \\
   h_n(s)  &= (1-p) \al G'\big(1-\al + \al g_{n-1}\big) h_{n-1}.
 \end{align}
Let
\be \label{e:vndef0} 
v_n = 
\begin{cases}
 (1-\al) \lam (\lam \al)^{n-1} & \hbox{ if $1 \le n \le b$}, \\
  (1-\al) \lam (\lam \al)^b h_{n-b-1} & \hbox{ if $n \ge b+1$}.
\end{cases} 
\ee

\begin{theorem} \label{T:Main1}
Let $p>0$. Then the process after contact tracing becomes extinct if and only if
\be \label{e:esum0}
\sum_{n=1}^\infty  v_n \le 1
\ee
\end{theorem}

We do not know of any case where we can give exact expressions for $(v_n)$,
but unless $p$ is small the series for $(h_n)$ converges very rapidly.

\sms
If $\al$ is not large enough to control the epidemic, we can still ask how quickly 
$Z^{CT}$ grows. Define the 
{\em Malthusian parameter}  of $(v_n)$  to be the unique $\th$ such that
\be  \label{e:Mal}
 \sum_{n=1}^\infty e^{-n \th} v_n = 1. 
\ee
We are concerned with the case when $\sum v_n>1$,  so $\th>0$.

\begin{theorem} \label{T:grow} 
On the event that $Z^{CT}$ survives, we have with probability 1,
\be \label{e:ZCT1}
 \lim_{n \to \infty} \frac{ \ln Z_n^{CT} }{n} = \theta.
\ee
\end{theorem}

An easy argument (see Lemma \ref{L:mon}) shows
that there is a function $e_b(p)$ such that the epidemic 
becomes extinct with probability one if $\al > e_b(p)$, and survives with positive probability
if $\al < e_b(p)$. (The function $e_b(p)$ depends on the offspring distribution of the original branching process.)
In Section \ref{S:crit} we study the function $e_b(p)$ close to the  critical points where it crosses
the axes. 

There is a very extensive applied epidemiological literature on contact tracing, but we did not
find very many papers containing exact calculations. One  quite closely related paper is
\cite{BKN}, which looks at forward contact tracing for a continuous time epidemic.
In the case of an exponential infection time they obtain an exact formula for the 
reproduction number of the epidemic after contact tracing.
As in our paper, they look at the discrete time process of untraced individuals --
called there unnamed individuals. \cite{BKN2} looks at some extensions of this model,
including allowing a latent period, and tracing delays.
The papers \cite{KFH, MK, MKD} also all consider various kinds of contact tracing 
for a continuous time branching processes.

\section{The contact tracing process }  \label{S:CT}

We need to keep track of not just the size of the original branching 
process, but also
its genealogical structure. So we write $\Lam_0 =\{0\}$, $\Lam_n = \bN^n$ and set
$$ \Lam = \bigcup_{n=0}^\infty \Lam_n, \q  \Lam[0,N] = \bigcup_{n=0}^N \Lam_n.$$
A point $x = (x_1, \dots ,x_n) \in \bN^n$ represents a potential individual in the $n$th generation; we write
$|x|=n$, say that $x$ is in generation $n$, and define its ancestor to be
$a(x)=(x_1, \dots x_{n-1})$. 
For $x \in \bN$ we set $a(x) =0$.
We take $\Lam$ to be a graph with edge set 
$$ E_\Lam = \Big\{ \{x,a(x)\}: x \in \Lam -\{0\} \Big \}. $$

Let $(\xi_x, x \in \Lam)$ be independent  r.v.~with distribution $(p_k)$.
We will define $\eta: \Lam \to \{0,1\}$, and set $\Gam = \{ x \in \Lam: \eta(x) =1\}$;
this will be the set of individuals in the original process.
We define $\eta(0)=1$. Once  $\eta$ is defined on $\Lam[0,n]$ we extend it to $\Lam[0,{n+1}]$ as follows. 
Let $x =(x', x_{n+1})\in \Lam_{n+1}$; we set $\eta(x)=1$ if $\eta(x') =1$ and $x_{n+1} \le \xi_{x'}$, and
take $\eta(x)=0$ otherwise.
Let $Z_n = |\{ x \in \Lam_n: \eta(x) =1 \}|$, so that $(Z_n)$ 
is a branching process with offspring distribution $(p_k)$.  
 We consider $\Gam$ as a graph with edge set
$E_\Gam = \big\{ \{x,a(x)\}: x \in \Gam, x \neq 0 \}$. 

Let  $b \in \bZ_+$,  and probabilities $p \in [0,1]$ and $\alpha \in [0,1]$.
Define i.i.d. random variable $\eta_D(x)$ with a {\tt Ber$(p)$} distribution,
and i.i.d. $\eta_T(x)$  with a {\tt Ber$(\alpha)$} distribution.
If $\eta_D(x) =1$ we say $x$ is {\em detectable} or {\em detected}.
If $\eta_T(x) =1$  then we say the edge $\{x,a(x)\}$ is  {\em open} or {\em traceable}.
Thus $(\eta_T)$ defines a bond percolation process on $\Gam$ -- see \cite{G}. 
A path (i.e. sequence of edges) in $(\Gam,E_\Gam)$ is traceable if each edge in the
path is traceable. For $x \in \Gam$ we write $\sC(x)$ for the set of $y \in \Gam$ such that
$x$ and $y$ are connected by an traceable path; we call
$\sC(x)$ the {\em connected cluster} of $x$ (at time $n$). We always have $x \in \sC(x)$.

The contact tracing procedure operates as follows.
At time $n\ge 0$ we will define a subset $A_n$ of $\Gam \cap \Lam[0,n]$.  
If $b>0$ or $b=0$ and $\eta_D(0)=0$ we take $A_0=\{0\}$. If $b=0$
and $\eta_D(0)=1$ then we set $A_0=\emptyset$. (In this case the founding
individual is detected at time 0 and the process immediately becomes extinct.)

To construct $A_n$ from $A_{n-1}$, let
$$ A^*_n = A_{n-1} \cup \{ x \in \Gam\cap \Lam_n: a(x) \in A_{n-1} \}; $$
thus $A^*_n$ is $A_{n-1}$ together with the  offspring of the individuals in $A_{n-1}\cap \Lam_{n-1}$.
We now look at the r.v. $\eta_D(y)$ for $y \in A^*_n \cap \Lam_{n-b}$, and 
if  $\eta_D(y)=1$ then we remove $\sC(y)$ from $A^*_n$. 
Thus we set 
$$ A^R_n = \bigcup\{ \sC(y)\cap \Lam_0^n: y \in A^*_n \cap \Lam_{n-b}, \eta_D(y)=1 \},
\q A_n = A^*_n - A^R_n. $$
The {\em current generation} of the process $(A_n)$ is $A^{\rm CG}_n = A_n \cap \Lam_n$.
The size of the current generation is
$$ Z^{CT}_n = |  A^{\rm CG}_n |. $$
Since $A^*_n \subset \Gam \cap \Lam_n$, $Z^{CT}_n \le Z_n$. 
We call the process $A=(A_n, n \ge 0)$ the $(b,p,\alpha)$--{\em contact tracing process}
or CTP$(b,p,\al)$. The parameter space is
\be \label{e:parsp}
 \sP = \{ (b,p,\al):  b \in \bZ_+, p \in [0,1], \al \in [0,1]\}. 
\ee 
We write $\bE_{b,p,\al}$ for expectations when we wish to emphasize the dependence on the parameters.

Note that not all points are removed -- for example if $\eta_D(y)=\eta_T(y)=0$ and $y$
has no descendants then $y \in A_n$ for all large $n$.
If however for some $n$ we have 
$A^{\rm CG}_n = \emptyset$ then $A^{\rm CG}_{n+k}=  \emptyset$ for all $k \ge 0$.

\begin{figure}[h!]
\includegraphics[width=\textwidth]{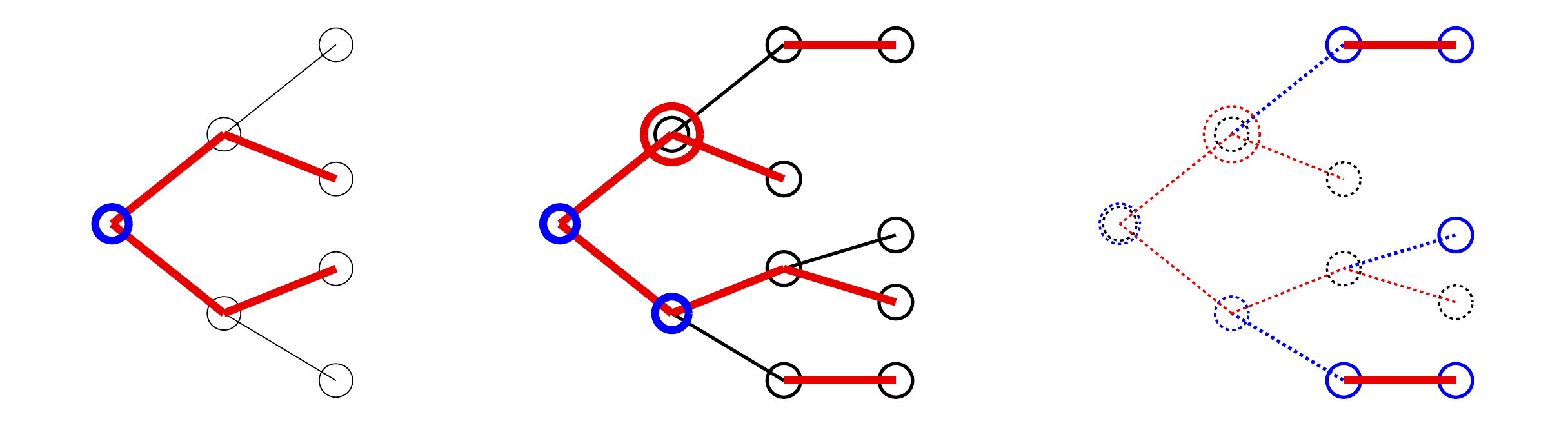}
\caption{The process $A_n$ at times $n=2$ (left), $n=3$ before tracing (center) and 
$n=3$ after tracing (right). Red (thick) lines are traceable edges.
At time 2 the root is not detected. At time 3 one vertex in generation 1 is detected (large red circle),
and it and the 5 vertices connected to it are removed. Hence one member of the current
generation is removed.} 
\end{figure}

To clarify our terminology we define what we mean by survival or extinction
of a random process.

\begin{definition}
{\rm Let $X=(X_n, n \ge 0)$ be a random process on $\bZ_+$ 
with the property that  $\{X_n =0\} \subset \{X_{n+1} =0\}$. 
We define 
$$ \{ X \hbox{ becomes extinct}\} = \bigcup_{n=0}^\infty \{X_n =0\}. $$
If $\bP(  X \hbox{ becomes extinct })=1$ we say {\em $X$ becomes extinct}, and if
$\bP(  X \hbox{ becomes extinct })<1$ we say {\em $X$ survives with positive probability}
or {\em survives wpp.} Sometimes we will shorten `survives wpp' to `survives'. 
For a set valued process such as $A^{\rm CG}=(A^{\rm CG}_n)$ we define extinction or survival as whether 
or not the process $|A^{\rm CG}_n|$ becomes extinct, or survives.
} \end{definition}

We are interested in characterizing the set
of parameter values such that $Z^{CT}$  becomes extinct, and define
the following subsets of $\sP$:
\begin{align*}
 \sE &= \{ (b,p,\alpha): \hbox{ CTP$(b,p,\alpha)$ becomes extinct} \} , \\
\sS &= \{ (b,p,\alpha): \hbox{ CTP$(b,p,\alpha)$ survives wpp} \}.
\end{align*}

\sms We begin with some easy properties of the sets $\sE$ and $\sS$. 

\begin{lemma} \label{L:easy}
(a) If $b=0$ and $\lam(1-p) \le 1$ then $Z^{CT}$ becomes extinct. \\
(b) If $b\ge 1$ and $\al=0$ then $Z^{CT}$ survives wpp. \\
(c) If $\lam(1-p)(1-\alpha)  > 1$  then $Z^{CT}$ survives wpp. \\
(d) If $p=0$ then  $Z^{CT}$ survives wpp. \\
(e) If $\alpha=1$ and $p>0$  then $Z^{CT}$ becomes extinct.

\end{lemma}

\proof (a) Let $(A'_n)$ be the process $A$, but where at time $n+1$ only the detected 
individuals are removed. Then $A_n \subset A'_n$, and the process $Z'_n=|A'_n \cap \Lam_n|$ 
is a simple branching
process with offspring distribution mean $(\lam (1-p))$; if $\lam(1-p) \le 1$
then $Z'$ becomes extinct, and therefore $Z^{CT}$ also becomes extinct. \\
(b) In this case detection never removes individuals in the current generation, so
$|A^{\rm CG}_n|$ is just a branching process with offspring distribution $(p_k)$, 
and therefore it survives wpp.\\
(c) If we just consider the new points $y$ which satisfy $\eta_T(y)=\eta_D(y)=0$, we have 
a process $A^{UN}$ which is smaller than $A$, and is a branching process with 
offspring distribution with mean $\lam (1-p)(1-\al)$, and  if this process survives then
$A$ survives. \\
(d) and (e) are clear. \qed

We will treat the case $p=1$ when $b\ge 1$ in Lemma \ref{L:border2} below.

\begin{remark}
{\rm The `heavy tailed' situation when $\sum p_k = \infty$ does not seem to be interesting
in this context. If $p=0$ or $\al=1$ then Lemma \ref{L:easy}(d),(e) still hold.
If $p=1$ and $\al<1$ the process becomes extinct if $b=0$ and survives otherwise, while if
$p \in (0,1)$ and $\al<1$ then the argument of Lemma \ref{L:easy}(c) implies that
the process survives wpp.
}\end{remark}

\begin{lemma} \label{L:mon} (Monotonicity).
If $(b,p,\alpha) \in \sE$ and $b'\le b$, $p' \ge p$ and 
$\alpha' \ge \alpha$ then $(b',p',\alpha') \in \sE$. 
\end{lemma}
 
\proof We write $\eta^{(p)}_D(x)$ and $\eta^{(\alpha)}_T(x)$ for the detection and tracing
processes with parameters $p$ and $\alpha$ respectively, and couple them so that they
are monotone in $p$ and $\alpha$ respectively. 
Write $A^{(b,p,\alpha)}$ for the associated contact tracing process. Then
the construction of $A$ gives that $A_n^{(b,p',\alpha')} \subset A_n^{(b,p,\alpha)}$.
Monotonicity in $b$ is also clear. \qed

Using Lemma \ref{L:mon}  we see there exists a function
$e_b: [0,1] \to [0,1]$ such that if $\al < e_b(p)$ then $(b, p, \al ) \in \sS$ and if 
$\al > e_b(p)$ then $(b, p, \al ) \in \sE$. Using Lemma \ref{L:easy} we see
that 
\be
  e_b(p) \ge 
\begin{cases}  
  1- \frac{1}{\lam(1-p)}, \, & \hbox{ if }  0\le p \le 1-\lam^{-1},\\
  0, & \hbox{ if }  1-\lam^{-1}\le p \le 1. \\
 \end{cases} 
\ee  
From Proposition \ref{P:e0-mglb} we will obtain the bound
\be  \label{e:e0mb-a}
 e_0(p) <  1 - p/\lam \q \hbox{ for } 0<p\le 1. 
 \ee

Lemma \ref{L:easy} covers the cases when $p=0$, so from now on we assume
\be \label{e:assume}
  p \in (0,1]. 
\ee

For $x \in \Gam \cap \Lam_n$  let $V_x$ be the number of offspring of $x$,
and $V^T_x$ and $V^U_x$  be the number of traceable and untraceable offspring.
So $V_x = V^T_x+V^U_x$ and
$$ \bE V^T_x = \lam \al , \q  \bE V^U_x = \lam(1-\al). $$
Write $p^T_k = \bP( V^T_x =k)$, and let 
$G_T(u)$ be the associated p.g.f. Then since the law of
$V^T_x$ conditional on $V_x$ is {\tt Binom$(V_x, \al)$},
\be  \label{e:GTdef}
 G_T(u) = \sum_{k=0}^\infty u^k p^T_k = G( (1-\al) + \al  u ), \\
\ee

We can decompose $A_n$ into a collection of disjoint connected 
traceable clusters $\sC^1, \dots , \sC^{k_n}$, and an analysis
of the evolution of a traceable cluster is a key step in our proofs.

We will call $y\in  \Gam$ with $\eta_T(y)=0$  a {\em cluster seed}.
The edge $\{a(y),y)\}$ is untraceable, and so if $|y|=n$ and $y \in A_n$ then $y$
will only be removed 
from the process $(A_{n+k}, k\ge 0)$ if some descendent of $y$
is detected and is connected to $y$ by a traceable path. 
A key observation  is that
the evolution of the part of $(A_{n+k}, k \ge 0)$ containing $y$ and its
descendants is independent of the rest of  $(A_{n+k}, k \ge 0)$.

We look at the evolution in time of a traceable cluster, and for simplicity
consider the cluster started at the root $0$. 
Initially we consider the growth of the cluster without any detection.
Let $\sV_0 = \{0 \}$, and $(\sV_k, k \ge 1)$ be the cluster
at subsequent times, given by
$$ \sV_n = \{ x \in \Gam \cap \Lam_n: a(x) \in \sV_{n-1}, \eta_T(x) = 1 \}. $$
Let $V^T_0=1$,  $Y^U_0=Y^U_0=0$, and  for $n \ge 1$ let
\begin{align}
 \label{e:VnU}
 &V_n^U = \sum_{x \in \sV_{n-1}} V^U_x ,   &Y^U_n = \sum_{k=1}^n V^U_k, \\
  &V_n^T = \sum_{x \in \sV_{n-1}} V^T_x= |\sV_n|,  &Y^T_n = \sum_{k=0}^n V^T_k. 
 \end{align} 

We now introduce detection for this cluster. Let
$S$ be the generation number of the first detected vertex in the cluster. So
$$ \{ S =n \} = \{ \eta_D(y) =0 \hbox{ for all } y \in \bigcup_{k=0}^{n-1}\sV_k \} \cap
  \{ \sum_{y \in \sV_n } \eta_D(y) \ge 1 \}. $$
Write $\wt V^T_n$ for the number  of points in the cluster
 taking detection into account, and $\wt V^U_n$ for the corresponding number of cluster seeds.
 We have
 \be \label{e:wtV}
   \wt V^T_n = V^T_n 1_{(n \le S+b-1 )}, \q \wt V^U_n = V^U_n 1_{(n \le S+b)}.
 \ee
 (The difference between the two expressions above is because the detection process forces
 $\wt V^T_{S+b}=0$, while the cluster may still produce cluster seeds in generation $S+b$.)
Thus the total number of cluster seeds produced by the cluster starting from 0 is
\be
\wt Y^U_\infty  =  \sum_{n=1}^\infty \wt V^U_n  = \sum_{n=1}^\infty V^U_n 1_{(n \le S+b)}.
\ee 
If  $b=0$ and $\eta_D(0)=1$ then $0$ is detected and $S=\wt Y^U_\infty=0$. 
For $x \in \cup A_n$ let  $\sY(x)$ be the set of cluster seeds produced by the traceable cluster
with cluster seed $x$. 
Set
\be \label{e:vndef}  
v_n =   \bE( \wt V^U_n ) \q \hbox{ for } n \ge 1.
\ee

\begin{thm} \label{T:EW1}
Let $p>0$. The process $Z^{CT}$ becomes extinct if and only if 
$$ y_b(p,\al) = \bE_{b,p,\al}( \wt Y^U_\infty )= \sum_{n=1}^\infty v_n  \le 1.$$
\end{thm}

 \proof 
Let $\sX_0=0$, $X_0= 1$, and define for $n \ge 1$
$$ \sX_n = \bigcup_{x \in \sX_{n-1} } \sY(x), \q X_n  = |\sX_n|. $$
Then $X=(X_n)$ is a branching process with offspring distribution equal in law  to $\wt Y^U_\infty$.
Since  $  \bP( \wt Y^U_\infty  =1)  <1$ the process $X$ 
becomes extinct with probability 1 if and only if $\bE X_1 \le 1$, i.e. if and only if $\bE \wt Y^U_\infty \le 1$. 

It remains to show that $A^{\rm CG}$ (a.s.) becomes extinct if and only if $X$ becomes extinct.
 If $A^{\rm CG}$ (a.s.) becomes extinct  then the total family size
$\sum_n | A^{\rm CG}_n|$ is finite, so the total number of cluster seeds is finite and thus
$X$ becomes extinct.  On the other hand, if $X$ becomes extinct then (as $p>0$)
each traceable cluster is finite, so the total family size
$\sum_n | A^{\rm CG}_n|$ is finite, and thus $A^{\rm CG}$ becomes extinct. \qed

\begin{proposition} \label{P:e0-mglb} 
If $b=0$ and $\al \ge 1 - p/\lam$ then the process $Z^{CT}$ becomes extinct. Hence
\be  \label{e:e0mb}
 e_0(p) \le 1 - p/\lam \q \hbox{ for } 0<p\le 1. 
 \ee
\end{proposition}

\proof 
We consider the traceable cluster $(\sV_n)$ started at the origin.
We explore the points in the traceable cluster, starting with the  origin,
and completing the exploration of each generation before starting on the next.
Let $X_1, X_2, \dots$ be the exploration process; 
the r.v. $X_i$ take values in the genealogical space $\Lam$.
This sequence is finite or infinite according to whether the cluster is finite or infinite.
The total family size of $(\sV_n)$ is $Y=Y^T_\infty$.  

We add a cemetery point $\pd$, define $\eta_D(\pd)= V^U_\pd=0$ and if $Y<\infty$ we set
$X_n = \pd$ for $n>Y$. Set
$$ \sF_n = \sigma(X_j,  \eta_D(X_j), V^U_{X_j}, 1\le j \le n ). $$
Let $T_E = \min\{n \ge 1: X_n = \pd\}$; then $T_E=Y+1$ and is a stopping time with respect
to $(\sF_n)$. Let $T_D = \min\{ n\ge 1: \eta_D(X_n) =1\}$, and $T = T_D \wedge T_E$.
Set $U_0=D_0=M_0=0$, and for $n \ge 1$ let
$$ U_n = \sum_{k=1}^n V^U_{X_k}, \q D_n = \sum_{k=1}^n \eta_D(X_j), 
\q M_n = p U_n - \lam(1-\al) D_n. $$
Then $M$ is a martingale with respect to $(\sF_n)$. Hence for any $n \ge 0$,
$$ 0 = M_0 = \bE( M_{T\wedge n}) = p \bE (U_{T\wedge n} )  -   \lam(1-\al) \bE( D_{T\wedge n}) . $$
The total number of cluster seeds produced by the traceable cluster
up to its extinction due to a detection is $\wt Y^U_\infty$. 
Note that $\wt Y^U_\infty$ need not include the all r.v. $V^U_x$ for $x$ in the final generation, 
and so $\wt Y^U_\infty \le U_T$.  We have $D_T \le 1$. Using Fatou's Lemma
$$  \bE(\wt Y^U_\infty) \le  \bE(U_T)  \le \lim_{n\to \infty}  \bE U_{T \wedge n} 
= \lam(1-\al) p^{-1}  \lim_{n\to \infty}   \bE D_{T \wedge n} \le  \lam(1-\al) p^{-1} . $$
The conclusion  is now immediate from Theorem \ref{T:EW1}. \qed

\section{ Calculation of the number of cluster seeds}  \label{S:calc}

Set
$$  \lam_T = \lam \al , \q  \lam_U = \lam(1-\al),$$
and let
\be \label{e:wndev}
   w_n = \bE (   V^T_n (1-p)^{Y^T_n} ), \q n \ge 0.
\ee

\sms
\begin{lemma}
We have 
\begin{align*}
 \bE_{b,p,\al} ( \wt V_n^U ) &=   \lam_U \lam_T^{n-1}   \q \hbox{ for } 1\le n \le b,\\
     \bE_{b,p,\al} ( \wt V_{n}^T  ) &=    \lam_T^n   \q \hbox{ for } 0\le n \le b-1, \\
  \bE_{b,p,\al} ( \wt V_n^U ) &=   \lam_U \lam_T^b w_{n-b-1}
   \q \hbox{ for } n \ge b+1, \\
      \bE_{b,p,\al} ( \wt V_{n}^T  ) &=    \lam_T^b w_{n-b}
       \q \hbox{ for } n \ge b.
\end{align*}
Hence
\be \label{e:ybdef}
y_b(p,\al)= \bE_{b,p,\al}(\wt Y^U_\infty) = 
\sum_{n=1}^b \lam_U \lam_T^{n-1} +  \lam_U \lam_T^b 
  \sum_{n=0}^\infty  w_n.  
\ee

\end{lemma}

\proof 
The first two expressions follow easily from \eqref{e:wtV}. 
Let $\sF^T_n = \sigma( V^T_j, 0 \le j \le n, \sV_i, 0\le i \le n )$.  We have
$ \{ S \ge n \} = \{ \eta_D(x)=0 \hbox { for all } x \in \cup_{i=0}^{n-1} \sV_i \}$.

Now let $ n \ge b+1$. 
The r.v. $1_{(S\ge n-b)}$ and $V^U_{n}$ are conditionally independent given
$\sF^T_{n-b-1}$. So
\begin{align*}
  \bE ( \wt V_n^U |  \sF^T_{n-1}   ) = 
   \bE (  V_{n+b}^U 1_{(n \le S)} | \sF^T_{n-1} ) 
  =   \bE (  V_{n+b}^U  | \sF^T_{n-1} ) \, \bE (   1_{(n \le S)} | \sF^T_{n-1} ). 
  \end{align*}
Then
$$  \bE (  V_{n+b}^U  | \sF^T_{n-1} ) =  \bE \big(  \bE( V_{n+b}^U | \sF^T_{n+b-1} )  | \sF^T_{n-1} \big) 
 = \bE \big(  \lam_U V_{n+b-1}^T  | \sF^T_{n-1} \big) = \lam_U \lam_T^b V^T_{n-1}, $$
while
$$ \bP( S \ge n | \sF^T_{n-1} ) = (1-p)^{V^T_0 + V^T_1+ \dots + V^T_{n-1} } = (1-p)^{Y^T_{n-1}}. $$
Combining these equalities and taking expectations gives the expression for $\bE ( \wt V_n^U )$.

Similarly if $n \ge b$ we have
\begin{align*}
  \bE ( \wt V_n^T |  \sF^T_{n-b}   ) &= 
   \bE (  V_{n}^T 1_{(S \ge n -b+1 )} | \sF^T_{n-b} )   \\
  &=   \bE (  V_{n}^T  | \sF^T_{n-b} ) \, \bE (   1_{(S \ge n-b+1)} | \sF^T_{n-b} )  
  = \lam_T^ b V_{n-b}^T (1-p)^{Y^T_n}.
  \end{align*}
\qed

To calculate $w_n$ set for $s, t \ge 0$, $n \ge 0$
\be
  H_n(s,t) = \bE ( s^{Y^T_n} t^{V^T_n} ). 
\ee 
Note that $H_0(s,t) =st$.
Then
$$  \frac{\pd} {\pd t}  H_n(s,t) = \bE ( V^T_n s^{Y^T_n} t^{V^T_n-1} ), $$
and so
\be \label{e:pdH}
   \bE (  V^T_{n} (1-p)^{Y^T_{n}} ) =  \frac{\pd} {\pd t}  H_n(1-p,1). 
 \ee  
Since $p>0$ and $V_n^T \le Y_n^T$, the series for $H_n(1-p,t)$ converges in a neighbourhood
of 1,  and there is no problem taking the derivative in \eqref{e:pdH}.
A standard first generation branching process decomposition  gives
\begin{align*}
  H_n(s,t) = s G_T (H_{n-1}(s,t)) = s G( 1-\al + \al H_{n-1}(s,t)).
 \end{align*}
Thus 

 \begin{align} \label{e:Hn}
\frac{\pd} {\pd t}  H_n(s,t) &= s \al G'(  1-\al + \al H_{n-1}(s,t) ) \frac{\pd} {\pd t}  H_{n-1}(s,t).
 \end{align}
 Setting
$$ g_n (s) = H_n(s,1), \q h_n(s) = \frac{\pd} {\pd t}  H_n(s,t), $$
we have the system of equations, for $s \in [0,1]$,
\begin{align}
   g_0(s) &=h_0(s)=s, \\
   g_n(s) &= s G\big( 1-\al + \al g_{n-1}(s)\big), \\
   h_n(s)  &= s \al G'\big(1-\al + \al g_{n-1}(s)\big) h_{n-1}(s),
 \end{align}
 and
 \be \label{e:wnval}
  w_n = h_n( 1-p). 
 \ee
 
\sm {\bf Proof of Theorem \ref{T:Main1}.} 
This follows from Theorem \ref{T:EW1},  \eqref{e:ybdef} and \eqref{e:wnval}. 
We also have that $v_n$ as defined by \eqref{e:vndef} satisfies  \eqref{e:vndef0}.
\qed
 
 Note that the functions  $g_n$ and $h_n$ depend on $\al$ but not on $p$.
We collect some properties of these functions. 

\begin{lemma} \label{L:basic}
(a) The functions $g_n(s)$, $h_n(s)$ are strictly increasing and continuous for $s\in [0,1]$.\\
(b) For $s\in [0,1)$, the sequence  $g_n(s)$ is strictly decreasing in $n$.
The limit $g_\infty(s)$ satisfies $g_\infty = s G_T(g_\infty)$ and 
$s G'_T(g_\infty(s)) < 1$. \\
(c) Let $c_1(p)= (- e \log(1-p))^{-1}$.   We have 
$$ h_n(1-p) \le c_1(p) (1-p)^{n} \hbox{ for } n \ge 1. $$
Hence for $p>0$
\be \label{e:hsumub}
  \sum_{n=0}^\infty h_n(1-p) \le \frac{ 1} { ep \log(1/(1-p)) }. 
\ee
(d) For each $b \in \bZ_+$ the function $y_b(p,\al)$ is continuous
in the region $(p,\al) \in (0,1] \times [0,1]$.
\end{lemma}

\proof (a) This is clear from the definition. \\
(b) Note that $G_T$ is strictly monotone.
Fix $s \in (0,1)$. Then $g_1 = sG_T(s)< s G_T(1) =s$. If $g_{n-1}<g_{n-2}$
then $g_n = s G_T(g_{n-1}) < s G_T(g_{n-2}) = g_{n-1}$. Thus $g_n$ is decreasing
in $n$, and as $G_T$ is continuous the limit must satisfy  $g_\infty = s G_T(g_\infty)$.

Let $f_1(x) =x$ and $f_2(x) =s G^T(x)$. Then
$0=f_1(0) < f_2(0) = sG_T(0)$, while $f_1(1) =1 > s = f_2(1)$. Thus
$f_1(x)=f_2(x)$ has a solution in $[0,1]$, and as $G_T$ is strictly
monotone it follows that this solution is unique, and therefore equals $g_\infty$.
By the mean value theorem there exists $\xi \in (g_\infty,1)$ such that
$s-g_\infty = f_2(1)-f_2(g_\infty) = (1-g_\infty) f'(\xi)$. As $f_2'$ is monotone
we have $s-g_\infty  \ge  (1-g_\infty) f'(g_\infty)$, and thus $f'_2(g_\infty)<1$. \\
(c) Note that if $V^T_n>0$ then $Y_n^T \ge n+1$, and that for
$x \ge 0$ we have $x (1-p)^x \le c_1$.  So
$$  h_n(1-p) = \bE V_n^T (1-p)^{V^T_n+ Y^T_{n-1}} \le c_1 (1-p)^{n}. $$
(d) If $K$ is  a compact subset of $(0,1] \times [0,1]$ then
by (c) the functions $h_n$ converge uniformly to 0 in $K$. It is straightforward to
verify that each term in the sum \eqref{e:ybdef} is continuous in $p,\al$; as the series
is uniformly convergent the limit is continuous. 
\qed

\begin{corollary}
We have $(b,p,e_b(p)) \in \sE$ if $p>0$.
\end{corollary}

\proof 
As $y_b$ is continuous in the region $p>0$, and extinction occurs if $y_b(p,\al)\le 1$,
extinction occurs in the critical case $\al = e_b(p)$. \qed

\begin{remark} 
{\rm  Note that $y_b$ is not continuous at $(0,1)$, since $y_b(0,1)=0$ while
for $\al \in (\lam^{-1},1)$ we have $y_b(0,\al)=\infty$.
Further $e_b(0)=1$, but $(b,0,1) \in \sS$, so the restriction to $p>0$ in the Corollary above
is necessary.
} \end{remark} 

The following Lemma handles the case $p=1$ and $b \ge 1$.

\begin{lemma} \label{L:border2}
(a) Let $b \ge 1$, and let
$$ \alpha_{\lam,b} = \inf \big\{ \al \ge 0: \lam (1-\al) \sum_{n=0}^{b-1} \lam^n \al^n \le 1 \big \}. $$
Then $\al_{\lam,b} \in (0,1)$, and for $p\in [0,1]$
CTP$(b,1,\al)$ survives wpp if $\al \in [0, \al_{\lam,b})$, and becomes extinct
if $\al \in [\al_{\lam,b},1]$.   \\
(b) Let $b \ge 1$. Then $e_b(p) \ge \al_{\lam,b}$ for $p \in [0,1]$.
\end{lemma}

\proof (a) Set 
\be \label{e:fbdef}
 f_b(\al) =  \lam(1-\al) \sum_{n=0}^{b-1} (\al \lam)^n. 
\ee 
We have $y_b(1,\al) = f_b(\al)$. Note that $f_b(0)=\lam$, $f_b(1)=0$, and $f_b$ is continuous.
Hence $\al_{\lam,b} \in (0,1)$, 
and the result follows using Theorem \ref{T:Main1} and Lemma \ref{L:mon}.\\
(b) This follows from (a) by monotonicity.
\qed

\begin{remark} 
{\rm In spite of the monotonicity given by Lemma \ref{L:mon}, 
the function $y_b(p,\al)$ is not monotone in $\al$. 
An easy  way to see this is to note that for $b \ge 1$ we have 
$$ \frac{\pd}{\pd \al} y_b(1,\al)\Big|_{\al=0}  = \lam^2 -\lam >0. $$
} \end{remark}

%\sms
%The following table gives $e_0(p)$ when the original branching process
%has various offspring distributions.
%
%\begin{table}[h!]
%  \begin{center}
%    %\caption{Your first table.}
%    \label{tab:table1}
%    \begin{tabular}{c c | c  c  c c c c c} 
%    % <-- Alignments: 1st column left, 2nd middle and 3rd right, with vertical lines in between
%%      \textbf{Value 1} & \textbf{Value 2} & \textbf{Value 3}\\
%Distribution & Mean      & $p=0.1$ & $p=0.2$ & $p=0.3$ & $p=0.4$ & $p=0.5$ & $p=0.6$ & $p=0.7$ \\
%      \hline
%%      \hline
%Poisson    &  2.5 & 0.9267 & 0.8417  &0.7401 & 0.6127 & 0.4366 & 0.00 & 0.00 \\ 
%Poisson    &  3.0 & 0.9377 & 0.8661  & 0.7816 & 0.6785 & 0.5446 & 0.3449 &0.00 \\
%Poisson     &  3.5 & 0.9447 &0.8815 &0.8075 &0.7181 &0.6050 &0.4476 &0.1516 \\
%Deterministic & 3  & 0.9414 &0.8745 &0.7960 &0.7002 &0.5747 &0.3809 &0.00 \\
%Geometric & 3.0 & 0.9162 &0.8280 &0.7321 &0.6222 &0.4852 &0.2891 &0.00 \\
%   \end{tabular}
%  \end{center}
% \caption{The function $e_0(p)$ for various offspring distributions. }
%\end{table}

\begin{figure}[h!]
\includegraphics[width=0.7\textwidth]{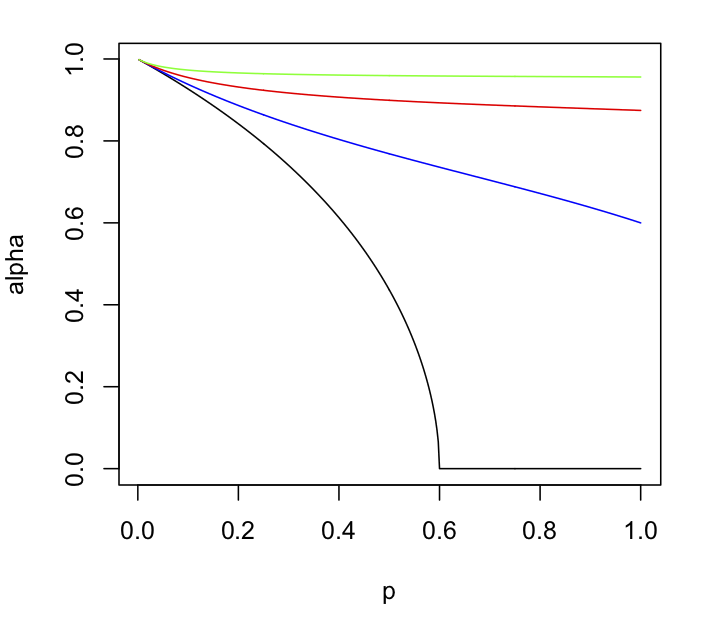}
\caption{ The functions $e_b(p)$ for $b=0,1,2,3$,  and a Poisson mean 2.5 offspring distribution.} 
\end{figure}

\section{Some limiting results}  \label{S:crit}

In the previous section we saw that we can obtain accurate numerical estimates on
$e_b(p)$ in regions of the parameter space when $p$ is bounded away from zero.
In this section we look at the some limiting behaviours of this function,
beginning with the easy case of $b \to \infty$.

\begin{proposition} There exists $b_0$, depending only on $\lam$, 
and a constant $c_3(p,\lam)$ such that for $b \ge b_0$,
\be
  c_3(p,\lam)\lam^{-b} \le 1-e_b(p) \le 2 \lam^{-b}.
\ee 
\end{proposition}

\proof 
By \eqref{e:esum0} and \eqref{e:hsumub} we have, writing $c_2(p)^{-1} = -ep \log(1-p)$, 
\be \label{e:w1}
f_b(\al) \le y_b(p,\al) \le f_b(\al) +  \lam(1-\al) (\al \lam)^b \sum_{n=0}^\infty w_n
\le f_b(\al) +  \lam(1-\al) (\al \lam)^b c_2(p). 
\ee
Set $\al_1 = \al_1(b) = 1 -2 \lam^{-b}$; then bounding from below the sum of the left side of \eqref{e:w1}
by its final term,
$$ y_b(\al_1,p) \ge \lam^b (1-\al_1)\al_1^b = 2 (1- 2\lam^{-b})^b > 2 (1 - 2b \lam^{-b})>1 $$
provided $4b < \lam^b$. Thus $(b,p,\al_1)\in \sS$ survives, and so by
monotonicity $(b,p,\al)\in \sS$ for all $\al \in [0,\al_1]$, and $e_b(p) \ge \al_1$.

For the lower bound we can assume that $b$  is large enough so that 
$\al_1 \lam-1 > \half (\lam-1)$. Then we obtain from \eqref{e:w1} that
$$ y_b(p,\al) \le \lam (1-\al) \lam^b \Big( 2 \lam(\lam-1)^{-1}  + c_2(p) \Big), $$
and taking $\al = 1 - c_3 \lam^{-b}$, with $c_3$ small enough, 
the  bound  follows. \qed

We now look at $e_0(p)$ close to the point $p_0 = 1-\lam^{-1}$;
for simplicity we restrict to the case when the offspring distribution has a finite second moment.

\begin{proposition} \label{P:2crit} 
(a) If 
\be \label{e:s1}
  \lam(1-p) (1-\al) \Big( 1 + \al (1-p) G'(1-\al p)  \Big) > 1
\ee
then $(Z^{CT}_n)$ survives. \\
(b)  Assume that $\sum k^2 p_k <\infty$. 
Then 
\be \label{e:p0lim}
  \lim_{ p \uparrow p_0} \frac{ e_0(p) } {(p_0-p)^{1/2} } = \frac{\lam}{p_0(1-p_0) G''(1)}. 
\ee
\end{proposition}

\proof 
(a) Note that 
$h_1(s) = s^2 \al G'( 1 -\al p)$. So
taking only the first two terms in the sum in \eqref{e:esum0} we have
$$ y_b(\al,p) \ge  (1-\al) \lam (1-p) \Big( 1 + (1-p) \al G'(1-\al p) \Big),$$
and the result is then immediate from Theorem \ref{T:EW1}. \\
(b) Write $t = p_0-p$ so that $\lam (1-p) = 1+ \lam t$.
Note that $G'(1-\al p) = \lam - \al p G''(1) + o(\al^2)$. 
We have
\begin{align*}
 \lam(1-\al) \sum_{n=0}^1 h_n(1-p)
 & = (1+\lam t) (1-\al) \Big( 1 + \al (1-p) G'(1-\al p) \Big)    \\
 & =  1+\lam t   -\al^2 -  \al^2 p_0(1-p_0) G''(1) + O(\al t) +   o(\al^2).     \
 \end{align*}
Further
\begin{align*}
  \lam(1-\al) h_2(1-p) &=  \lam(1-\al) \al^2 (1-p_0)^3 G'(1-\al p) G'( 1- \al + g_1(1-p) \al) \\
   &= \al^2 + O(\al^3, \al t ). 
  \end{align*}   
For $n \ge 3$ we have
$h_n(1-p) \le \al (1-p) \lam h_{n-1}(1-p)$, and thus we have
$$ y_0(p,\al) = 1 +\lam (p_0-p) -  \al^2 p_0(1-p_0) G''(1) + o(\al^2) + O(\al t) $$
 which gives the limit \eqref{e:p0lim}. \qed 

\begin{proposition}
Assume that $\sum k^2 p_k <\infty$. There exists $c_1, c_2$ such that
$$ e_0(p) \ge 1 - c_1 p \q  \hbox{ for } 0<p < c_2. $$
\end{proposition} 

\proof
Let $\al_0\in (0,1)$ be such that $\lam \al_0 >1$. Then it is easy to verify that
there exist $c_3, c_4$ such that for all $n \ge 1$ and $\al \in [\al_0,1]$,
$$ \bE V^T_n \ge c_3 \lam_T^n, \q \bE (V^T_n)^2  \le c_4 \lam_T^{2n}, 
\q \bE( Y_n^T) \le c_4 \lam_T^n.   $$
Write $e^{-p_1} = (1-p)$, and choose $c_2$ so that
$2p\ge p_1 \ge p$ for $p \in [0,c_2]$. Let $\al \in [\al_0,1]$. 
Choose $n$ so that $\lam_T^{n-1} < p^{-1} \le \lam_T^n$.
Then
$$ w_0(p,\al) \ge \lam_U  \bE V^n_T  (1-p)^{Y_n^T} . $$
By the second moment method
$$ \bP( V^T_n \ge \half \lam_T^n)  \ge \frac{ (\bE V^T_n )^2} { 4 \bE (V^T_n)^2 }
 \ge  \frac{ c_3^2 } {4 c_4}.   $$
Let $r = c_3^2/4c_4$. Then
$$ \bP( Y^T_n > 2 c_2 \lam_T^n/r ) \le \half r,$$
and thus 
$$ \bP( V^T_n \ge \half \lam_T^n ,  Y^T_n \le 2c_2 \lam_T^n/r  ) \ge \half r. $$
Hence writing $s=\half$, $t = 2c_2 /r$, 
\begin{align*}
 \bE V^n_T  (1-p)^{Y_n^T}  
 &\ge  \bE \Big( V^n_T  (1-p)^{Y_n^T} ;  V^T_n \ge s \lam_T^n ,  Y^T_n \le t\lam_T^n) \\
 &\ge \half r  s \lam_T^n \exp( { -p_1 t \lam_T^n}     )   \ge c p^{-1} e^{- c' \lam}.  
 \end{align*}
 Hence 
 $ w_0(\al,p) \ge  c_\lam (1-\al) /p, $ which gives that  $e_0(p) > 1 - c_1 p$. \qed
 
\section{Rate of growth of the process after contact tracing } 

In the case when $\al$ is not large enough to make $(Z^{CT}_n)$ extinct, it is of interest
to ask how quickly this process grows. 
We write
$$ Z^{CT}_n = |A^{\rm CG}_n|= Z^U_n + Z^T_n; $$
here $Z^U_n$ and $Z^T_n$ are the number of points $x \in A^{\rm CG}_n$ with
$\eta_T(x) =0$ and $\eta_T(x)=1$ respectively.
Recall from \eqref{e:vndef} the definition of $v_n$, and \eqref{e:Mal} that of the 
Malthusian parameter $\th$ of $(v_n)$.
We will be concerned with the case when $(Z^{CT}_n)$ survives wpp, so  $\th>0$. 

Let $(X_n)$ be the cluster seed process defined in the proof of Theorem \ref{T:EW1}.
Let $\sC_n$ be the set of cluster seeds at time $n$. Define a process
$R_n$ taking values in $\bZ_+^{\bZ_+}$ as follows. For a sequence  
$x \in \bZ_+^{\bZ_+}$ we write the $i$th component as $x(i)$. 
Set $R_0(0) = 1$, and $R_0(i) =0$ for $i \ge 1$.
Let $\sC_n$ be the set of cluster seeds at generation $n$, so
$$ \sC_n = \{ x \in A^{\rm CG}_n : \eta_T(x) =0 \}. $$
If $x$ is a cluster seed and $|x|=n$ write $\wt V^U_k(x) $ for the number of cluster seeds 
in generation $n+k$ produced by the traceable cluster with seed $x$.
Define
$$ R_n( k) = \sum_{m=0}^{n-1} \sum_{x \in \sC_m} \wt V^U_{n-m}(x) . $$
We have
$$ R_n(k) = R_{n-1}(k+1) + \sum_{x \in \sC_{n-1}} \wt V^U_{k+1}. $$
The process $R_n$ is a branching process, with infinitely many types.
For each $k \ge 1$ an individual of type $k$ produces with probability 1 one descendent
of type $k-1$. An individual of type 0 produces a random number of offspring,
with distribution $(\wt V^U_1, \dots)$. Let
$M_{ij}$ denote the mean number of offspring of a type $i$ individual. Then
\begin{align}
  M_{0,k} &= v_{k+1},  \q n \ge 0, \\
  M_{j,k} &= \delta_{k,j-1} ,\q  j \ge 1, \, k \ge 0. 
\end{align} 
It is straightforward to verify that if $u (n) = e^{-n \th}$, for $n \ge 0$ then $u$ is a right
eigenvector for $M$ and $M u = e^\th u$.
Thus we have 
$$ \bE( R_{n+1}(j) | R_n ) = \sum_{k=0}^\infty R_n(k) M_{k j}, $$
and
$$ Y_n = e^{-\th n} \sum_{k=0}^\infty R_n(k) u(k) $$
is a martingale with respect to $\sF^R_n = \sigma(R_0, R_1, \dots , R_n)$.
The following Lemma follows easily from Markov's inequality and the convergence of $Y$.

\begin{lemma} (a) We have
\be
 \bP( R_n(0) e^{-n \th}  > t) \le t^{-1}.
\ee
(b) With probability 1,
\be \label{e:lsZU}
 \limsup_{n \to \infty} \frac{ \ln R_n(0)}{n} \le \theta. 
\ee 
\end{lemma} 
 
Let $(X_n)$ be the cluster seed process, and let $F$ be the event that $(X_n)$ survives;
set $\bP( F) = p_F>0$.

\begin{lemma}
On the event $F$, a.s.
\be \label{e:liZU}
 \liminf_{n \to \infty} \frac{ \ln R_n(0)}{n} \ge \theta. 
\ee 
\end{lemma} 

\proof Let $\eps>0$.  
Let $N\ge 1$  and set
$$  v_n^{(N)}=  1_{(n < N)} \bE ( \wt V^U_n \wedge N), $$
and let $\th_N$ be the  Malthusian parameter for $(v_n^{(N)})$. 
As $v_n$ are bounded, we have for any $r>0$ that
$$ \lim_{N \to \infty}  \sum_{n=1}^\infty e^{-r n } v_n ^{(N)} =  \sum_{n=1}^\infty e^{-r n } v_n , $$
and it follows that we can choose $N$ large enough so that  $ \th-\eps < \th_N \le \th$.

Let $R^{(N)}$ be the process $R$, but with 
the offspring distribution for a type 0 individual given by  
$(N \wedge \wt V^U_1, \dots, N \wedge \wt V^U_N, 0, \dots)$.
We can do the truncation a.s. on the probability space, 
so that $R_n(k) \ge R_n^{(N)}(k)$ for all $n\ge 0$, $k \ge 0$.
As all individuals for  $R^{(N)}$ are of type $0,1, \dots , N-1$,
we can regard $R^{(N)}$ as a multi-type branching process with finitely many types.
Let $M^{(N)}$ be the matrix of means of the  truncated process $R^{(N)}$, 
and $u^{(N)}(k) = e^{-\th_N k}$ for
$0\le k \le N-1$. Then  $M^{(N)} u^{(N)} = e^{\th_N} u^{(N)}$. 
A classic limit theorem (see \cite[p. 192]{AN}) implies that 
$e^{- k \th_N } R^{(N)}_k$ converges a.s. and in $L^2$ to  $u_L W_N$, where
$u_L$ is the left eigenvector of $M^{(N)}$ with eigenvalue $e^{-\th_N}$, and $\bP(W_N>0)>0$. 
Thus a.s.~on $\{ W_N>0 \}$,
$$    \lim_{n \to \infty} \frac{ \ln R^{(N)}_n(0)}{n} \ge \theta_N > \theta -\eps.$$
Let 
\be \label{e:Geps}
G_\eps = \Big\{   \lim_{n \to \infty} \frac{ \ln R_n(0)}{n} \ge  \theta -\eps \Big\}. 
\ee
By the above, we have $\bP( G_\eps) = q>0$.

Let $m \ge 1$, and $T_m = \min \{ k \ge 0:  X_k \ge m \}$. 
We have $\bP( T_m <\infty | F) =1$.
Write $x_i, 1\le i \le m$ for the first $m$ cluster seeds in generation $T_m$, 
$R^i$ for the associated multitype branching process, and $G_\eps(i)$ for the event
defined by \eqref{e:Geps} for $R^i$. Then
$\bP( \cup G_\eps(i) | T_m <\infty) = 1 - q^m$. 
On the event $\cup G_\eps(i)$ we have
$\lim_{n \to \infty} \frac{ \ln R_n(0)}{n} \ge  \theta -\eps$, which implies that 
$\bP( G_\eps) \ge \bP(F) - q^m$.

Consequently, for any $\eps >0$ we have $\bP( G_\eps) \ge \bP(F)$, 
and \eqref{e:liZU} follows. \qed

\begin{remark}
{\rm Suppose in addition that $\sum k^2 p_k <\infty$. Let
$$ u_R(k) = e^-{\th} \sum_{k=n}^\infty e^{ -(k-n) \th } v_n. $$
Then $u_R$ is a left eigenvector for $M$ with eigenvalue $e^\th$, and $\sum_k u_R(k) u_L(k) <\infty$. 
By Criteriion III of \cite{VJ} the matrix $M$ satisfies Case II of \cite{M},
and it follows that
$$ \lim_n e^{-\th n} R_n \to v_R W \hbox{ in $L^2$}, $$
where $W$ is a random variable with $\bE(W)=1$ and $\bP(W>0)=p_F$.
} \end{remark}

\begin{theorem}
Let $(X_n)$ be the cluster seed process, and let $F$ be the event that $(X_n)$ survives.
On $F$ we have, a.s. 
\be \label{e:ZCTb}
 \lim_{n \to \infty} \frac{ \ln Z_n^{CT} }{n} = \theta.
\ee
\end{theorem}

\proof We recall the construction of $A^*_n$ from Section \ref{S:CT}. 
We consider first the case $b=0$. Let $M_n = |A^*_n \cap \Lam_n|$.
Note that
\begin{align*}
 R_n(0) &= \{ x \in A^*_n \cap \Lam_n: \eta_T(x) =0 \}, \\
 Z^U_n &= \{ x \in A^*_n \cap \Lam_n: \eta_T(x) =\eta_D(x) 0 \},
\end{align*}
and that $Z^{CT}_n \le M_n$.

Conditional on $M_n$ we have that $R_n \sim {\tt Bin}(M_n, 1-\al)$.  
Hence (see \cite{J}),
$$ \bP( | R_n - (1-\al) M_n | > t | M_n) \le 2 \exp( - t^2/3M_n ). $$

Let $\eps>0$, and set $m_n = e^{ n(\th + \eps)}$, $a = 2/(1-\al)$.
Then if $n \ge c_0(\al,p,\eps)$,
\begin{align*}
 \bP( M_n \ge 2 (1-\al)^{-1} m_n ) 
 &\le  \bP( M_n \ge 2 (1-\al)^{-1} m_n ), R_n < m_n) + \bP( R_n \ge m_n) \\
 &\le \bP( |(1-\al) M_n - R_n | > m_n , M_n \ge 2 (1-\al)^{-1} m_n) + e^{-\eps n} \\
&\le \exp( - m_n^2/ 6 (1-\al)^{-1} m_n ) + e^{-\eps n} 
\le 2  e^{-\eps n} .
\end{align*}
Thus we have
$$ \limsup_{n \to \infty} \frac {\ln M_n}{n} \le \th, \q a.s., $$
and the upper bound for $\ln Z_n^{CT}/n$ follows immediately.

As $M_n \ge R_n(0)$, it is immediate that on $F$
$$ \liminf_{n \to \infty} \frac {\ln M_n}{n} \ge \th, \q a.s.$$
The estimate 
$$ \bP( | Z^U_n - (1-\al)(1-p) M_n | > t | M_n) \le 2 \exp( - t^2/3M_n )$$
and a similar argument to that above gives that
$$ \liminf_{n \to \infty} \frac {\ln M_n}{n} \ge \th, \q \hbox{ a.s. on $F$}, $$
and as $Z^{CT}_n \ge Z^U_n$, the lower bound in \eqref{e:ZCTb} follows. \qed

\begin{figure}[!ht]
\includegraphics[width=0.7\textwidth]{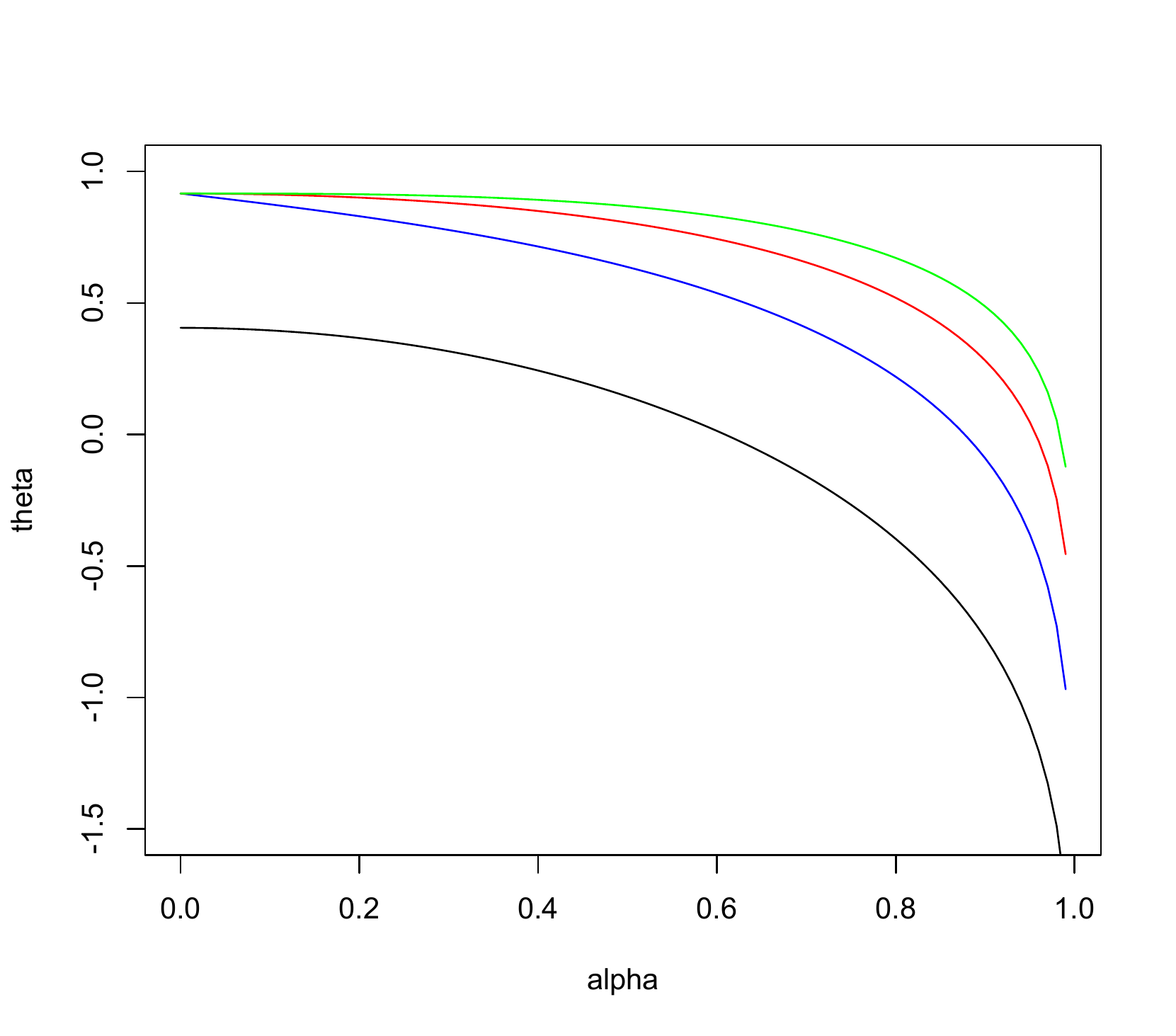}
\caption{ The Malthusian function $\th$ as a function of $\al$,
for a Poisson mean 2.5 offspring distribution, and $p=0.4$}
\end{figure}

\sm {\bf Question.} Let $(p_k)$ and $(p_k')$ be offspring distributions, 
and $\sE$ and $\sE'$ be the corresponding parameter sets for extinction.
Suppose that $(p_k')$ stochastically dominates $(p_k)$. It follows
that the branching process $Z'$ stochastically dominates $Z$. Is it the case that
$\sE' \subset \sE$?

\sm {\bf Acknowledgment.} 
This work was stimulated by a seminar in the B.C. COVID-19 group by Nils Bruin on contact tracing
-- see \cite{BCW} for the related paper.

\noindent MB: Department of Mathematics,
University of British Columbia,
Vancouver, BC V6T 1Z2, Canada. \\
barlow@math.ubc.ca


\begin{thebibliography}{}

\bibitem{AN} K.B. Athreya, P. Ney. {\em Branching processes.} Springer. 1972.

\bibitem{BKN}
F. G. Ball, E S. Knock, P. D. O'Neill.
Threshold behaviour of emerging epidemics featuring contact tracing.
{\em Adv. Applied Prob.} {\bf 43}, No. 4, (2011) pp. 1048--1065.

\bibitem{BKN2}
F. G. Ball, E S. Knock, P. D. O'Neill.
Stochastic epidemic models featuring contact tracing with delays. 
{\em Math. Biosciences} {\bf 266} 23--35 (2015) 

\bibitem{BCW} N. Bruin, C. Colijn, A. van der Waall.
Modelling the effects of epidemiological contact tracing, (in preparation).

\bibitem{G} G.R. Grimmett. {\it Percolation.} (2nd edition). Springer 1999.

\bibitem{H} T.E. Harris.  {Branching Processes.}
{\em Ann. Math. Statist.} {\bf19}, Number 4 (1948), 474--494.

\bibitem{J}  S. Janson
Large deviation inequalities for sums of indicator variables.
arXiv:1609.00533. (2016)

\bibitem{KFH} 
D. Klinkenberg,  C. Fraser,  H. Heesterbeek.
The Effectiveness of Contact Tracing in Emerging Epidemics.
PLoS One. 2006; 1(1): e12.
doi: 10.1371/journal.pone.0000012

\bibitem{M} S.-T. C. Moy.
Extensions of a limit theorem of Everett, Ulam and Harris on multitype branching 
processes to a branching process with countably many types.
{\em  Ann. Math. Statist.} {\bf 38}, Number 4 (1967), 992--999.

\bibitem{MKD} J. M\"uller, M.  Kretzschmar, K. Dietz.
Contact tracing in stochastic and deterministic epidemic models.
{\em Math. Biosciences} {\bf 164}, Issue 1, March 2000, Pages 39--64.

\bibitem{MK}  J. M\"uller, B. Koopmanna.
The effect of delay on contact tracing. 
{\em Math. Biosciences} {\bf 282} (2016) 204--214.

\bibitem{VJ}  D. Vere-Jones.
Ergodic properties of nonnegative matrices. I.
{\em Pacific J. Math.} {\bf 22}, (1967), 361--386.

\end{thebibliography}
\end{document}